\documentclass[11pt,a4paper]{article}
\usepackage{amsmath,amssymb}
\usepackage[title]{appendix}
\usepackage{relsize}
\usepackage{tensor}

\newcommand{\dis}{\displaystyle} 
\newcommand{\dop}{\mathrm{d}} 
\newcommand{\vsp}{\vspace{5mm}} 

\begin{document}
\title{\Huge \sf \bf An Algorithm for Calculating Terms of a Number Sequence using an Auxiliary Sequence}
\author{\sf Bengt M{\aa}nsson\\
\sf Kunskapsgymnasiet\\
\sf Gothenburg, Sweden}
\date{\sf\today}
\maketitle

\begin{abstract}
   \noindent Number sequences defined by a linear recursion relation are studied by means of generating functions. Indices of the terms in the recursion relation have arbitrary differences. In addition to formulas for the  $n$th term an algorithm is derived for calculating the $n$th term even without an expression in closed form.\footnote{This is version 2. An essential addition has been made at the beginning of page 3.}
\end{abstract}

The numbers are denoted $a_n, n=0,1,2,...$ and are defined by initial values

\[ \tensor{a}{_0}=\tensor{\alpha}{_0},\,\tensor{a}{_1}=\tensor{\alpha}{_1},\,\dots\tensor{a}{_{N-1}}
=\tensor{\alpha}{_{N-1}}\mbox{ för } N\geq 2 \]

and a recursion relation

\[ a_n=a_{n-1}+a_{n-N}\mbox{ for } n=N,\,N+1,\dots\ . \]

 The number $N$ has great influence on the complexity of derivations and calculations since it gives the degree of the polynomial in the denominator of the generating function.\\

\section{Generating function}

For generalities on using generating functions, see \cite{graham}, p 337 f and \cite{grimaldi}, p 445 f.\\

 A generating function for the number sequence $\left(a_n\right)_{n=0}^\infty$ is $f(x)=\sum_{n=0}^\infty a_n x^n$. Manipulation of indices and use of the recursion relation and the initial values yield

\begin{eqnarray*}
   f(x) &=& \sum_{n=0}^\infty a_{n+N} x^{n+N}+\sum_{n=0}^{N-1}\alpha_n x^n \\
   &=& \sum_{n=0}^\infty(a_{n+N-1}+a_n)x^{n+N}+\sum_{n=0}^{N-1}\alpha_n x^n \\
   &=& x\sum_{n=0}^\infty a_{n+N-1}x^{n+N-1}+x^N \sum_{n=0}^\infty a_n x^n+\sum_{n=0}^{N-1}\alpha_n x^n \\
   &=& x\sum_{n=0}^\infty a_n x^n-x\sum_{n=0}^{N-2}a_n x^n+x^N \sum_{n=0}^\infty a_n x^n+\sum_{n=0}^{N-1}\alpha_n x^n \\
   &=& (x+x^N)f(x)-\sum_{n=1}^{N-1}\alpha_{n-1} x^n+\sum_{n=1}^{N-1}\alpha_n x^n+\tensor{\alpha}{_0}
\end{eqnarray*}

from which we get

\begin{equation}\label{fx}
  f(x)=\frac{\dis \sum_{n=0}^{N-1}(\alpha_n-\alpha_{n-1}) x^n}{1-x-x^N}\ .
\end{equation}

Here we have defined $\tensor{\alpha}{_{-1}}=0$ in order to be able to write sums more concisely.\\

Since $a_n=f^{(n)}(0)/n!$ (taylor expansion) this gives

\begin{equation}\label{an1}
  a_n=\sum_{\ell=0}^{N-1}(\alpha_{\ell}-\alpha_{\ell-1}) C_{n\ell}
\end{equation}

where

\begin{equation}\label{dnl}
  C_{n\ell}=\dis\frac{1}{n!}\dis\frac{\dop^n}{\dop x^n}\dis\frac{x^{\ell}}{1-x-x^N}\Big|_{x=0}\ .
\end{equation}

\vsp
It will be seen later (section 3) that the right side expression depends on the differense $n-\ell$ only so $C_{n\ell}=D_{n-\ell}$ where $D_n$ is a certain function of $n$.

\vsp
Now, let $\delta$ denote the positive, real zero\footnote{For odd $N$ this is \emph{the only} real zero and it is always in the interval $]0,\,1[$.} of the polynomial $x^N+x-1$; $0<\delta<1$.\\

We now define another number sequence $(\kappa^N_{n,\ell})_{\ell=0}^{N-1}$ by

\[ (\delta^{N-1}+1)^n=\sum_{\ell=0}^{N-1}\kappa^N_{n,\ell}\delta^{\ell}\ \]

 where $N$ and $n$ are parameters specifying the number sequence and $\ell$ index, numbering the terms of the sequence.\\

The exponents of $\delta$ in the expansion of the left hand side can be reduced to at most $N-1$ since ${\delta}^N+\delta-1=0$, i e ${\delta}^N=1-\delta$. The powers of $\delta$ with indices $0,\,1,\,\dots\,N-1$ are linearly independent as follows from the fact that the polynomial $x^N-x-1$ is irreducible over the integers \cite{irred}, implying that this is the minimal polynomial of $\delta$. This motivates the identification of coefficients in the following.

\begin{eqnarray*}
  (\delta^{N-1}+1)^{n+1} &=& \sum_{\ell=0}^{N-1}\kappa^N_{n+1,\ell}\delta^{\ell}\ \\
   &=& \sum_{\ell=0}^{N-1}\kappa^N_{n,\ell}\delta^{N-1+\ell}+\sum_{\ell=0}^{N-1}\kappa^N_{n,\ell}\delta^{\ell} \\
   &=& \kappa^N_{n,0}\delta^{N-1}+\kappa^N_{n,1}+\kappa^N_{n,0}+\sum_{\ell=1}^{N-2}\kappa^N_{n,\ell+1}\delta^{\ell}
\end{eqnarray*}

giving

\begin{eqnarray*}
  \kappa^N_{n+1,N-1} &=& \kappa^N_{n,0} \\
  \kappa^N_{n+1,\ell} &=& \kappa^N_{n,\ell+1}\ (\ell=1,\,2,\dots,\,N-2\,;\,N\geq 3) \\
  \kappa^N_{n+1,0} &=& \kappa^N_{n,0}+\kappa^N_{n,1}
\end{eqnarray*}

and then

\begin{eqnarray*}
  \kappa^N_{n+1,N-1} &=& \kappa^N_{n,0}=\kappa^N_{n+1,0}-\kappa^N_{n,1} \\
   &=& \kappa^N_{n+2,N-1}-\kappa^N_{n-1,2} \\
   &=& \kappa^N_{n+2,N-1}-\kappa^N_{n-2,3}=\dots =\kappa^N_{n+2,N-1}-\kappa^N_{n-N+2,N-1}
\end{eqnarray*}

which finally gives

\begin{equation}\label{rekkappa}
   \kappa_n=\kappa_{n-1}+\kappa_{n-N}\quad .
\end{equation}

\vsp
Here we have excluded upper index and the second lower index since these have the same value in each term in each separate case (i e for each polynomial $x^N+x-1$ and each exponent $n$). Hence we have defined a new number sequence $(\kappa_n)_{n=0}^{\infty}$ where now $n$ is an index, numbering the terms of the sequence. It is easily seen that \eqref{rekkappa} holds also for $N=2$.

\vsp
Since \eqref{rekkappa} is the same recursion relation as the one for $a_n$ we get by means of \eqref{an1}

\begin{equation}\label{kappa}
  \kappa_n = \tensor{\kappa}{_0}D_n+\sum_{\ell=1}^{N-1}(\kappa_{\ell}-\kappa_{\ell-1})D_{n-\ell}
\end{equation}

and

\begin{equation}\label{an}
  a_n = \tensor{\alpha}{_0}D_n+\sum_{\ell=1}^{N-1}(\alpha_{\ell}-\alpha_{\ell-1})D_{n-\ell}\ .
\end{equation}

Independently of number sequence $a_n$ \footnote{The proof of this is somewhat complicated and is given in an appendix. In specific cases with small $n$ it is easily shown. For instance $(\delta^{N-1}+1)^2=\delta^{2N-2}+2\delta^{N-1}+1=\delta^N\delta^{N-2}+2\delta^{N-1}+1
=(1-\delta)\delta^{N-2}+2\delta^{N-1}+1=\delta^{N-1}+\delta^{N-2}+1$ so $\kappa_2=1$ and obviousily $\kappa_0=0$ and $\kappa_1=1$ (for all $N$).}

\[ \kappa_0=0,\ \kappa_n=1\ \mbox{för}\ 0<n\leq N-1\ \]

($N\geq 2$) which implies that \eqref{kappa} can be written

\begin{equation}\label{kappa2}
  \kappa_n =\sum_{\ell=1}^{N-1}(\kappa_{\ell}-\kappa_{\ell-1})D_{n-\ell}=D_{n-1}\ .
\end{equation}

Finally, substituting \eqref{kappa2} into \eqref{an} gives

\begin{equation}\label{ankappa}
  a_n = \tensor{\alpha}{_0}\kappa_{n+1}+\sum_{\ell=1}^{N-1}(\alpha_{\ell}-\alpha_{\ell-1})\kappa_{n+1-\ell}\ .
\end{equation}

In specific cases (given $N,\,\alpha_{\ell}$) this gives a relation between $a_n$ and $\kappa_n=\kappa^N_{n,N-1}$ by which $a_n$ can be calculated by expanding $(\delta^{N-1}+1)^n$ and noticing the coefficient for $\delta^{N-1}$, i e $\kappa^N_{n,N-1}$. Although this does not mean that we have an expression for $a_n$ in closed form, for specific, even large, values of $n$ the expansion often gives an easier calculation of $a_n$ than using a closed formula.\\

Notice that the exact expression for $D_n$ is not needed although, of course, it gives $a_n$ as an explicit function of $n$ by \eqref{an}.

\section{Examples}

Some useful relations to reduce polynomials in $\delta$ to polynomials with exponents in the interval $[0,\,N-1]$ and to increase/decrease the exponent $n$ in the expression of the form $(\delta^{N-1}+1)^n$:

\[ \delta^N=1-\delta,\ \frac{1}{\delta}=\delta^{N-1}+1\ . \]

\vsp
\textbf{Ex 1:} $N=2\mbox{ and } \alpha=0,\,1$ giving Fibonacci's number sequence.\\

In this case $\delta^2+\delta-1=0$, i e $\delta^2=1-\delta$ and $1/\delta=\delta+1$. The terms of the sequence are denoted $F_n$ and equation \eqref{ankappa} gives $F_n=\kappa_n$. We will find $F_{40}$.

\[ (\delta+1)^2=1-\delta+2\delta+1=\delta+2\ , \]
\[ (\delta+1)^4=(\delta+2)^2=1-\delta+4\delta+4=3\delta+5 \]

and similarly for $(\delta+1)^8$, $(\delta+1)^{16}$ och $(\delta+1)^{32}$. By means of this we get

\begin{eqnarray*}
  (\delta+1)^{40} &=& (\delta+1)^{32}\cdot(\delta+1)^8=(2178309\delta+3524578)(21\delta+34) \\
  &=& 45744489(1-\delta)+(74062506+74016138)\delta+3524578\cdot 34 \\
  &=& 102334155\delta+(\dots)\ .
\end{eqnarray*}

Thus, $F_{40}=102334155$.

\vsp
\textbf{Ex 2:} $N=3\mbox{ and } \alpha=0,\,1,\,2$\, .\\

In this case $\delta^3+\delta-1=0$, i e $\delta^3=1-\delta$, $\delta^4=\delta-\delta^2$ and $1/\delta=\delta^2+1$. Equation \eqref{ankappa} gives $a_n=\kappa_n+\kappa_{n-1}$. We will find $a_{16}$.

\[ (\delta^2+1)^2=\delta-\delta^2+2\delta^2+1=\delta^2+\delta+1\ , \]
\[ (\delta^2+1)^4=(\delta^2+\delta+1)^2
=\delta-\delta^2+\delta^2+1+2-2\delta+2\delta^2+2\delta=2\delta^2+\delta+3 \]

and similarly for $(\delta^2+1)^8$ and $(\delta^2+1)^{16}$. This gives

\[ (\delta^2+1)^{15}=(\delta^2+1)^{16}\cdot\delta=(189\delta^2+129\delta+277)\delta=129\delta^2+(\dots) \]

and, thus, we get

\[ a_{16}=189+129=318\ . \]

\vsp
\textbf{Ex 3:} $N=3\mbox{ and } \alpha=1,\,2,\,3$ and the recursion relation $a_{n+3}=a_{n+2}\cdot a_n$.

The number sequence\footnote{This example results from a discussion of which the next number should be in a sequence starting with 1, 2, 3, 3. The idea was that you can always find a ''rule'' giving any next number, for instance by adaption of a polynomial. One of my former collegues at Kunskapsgymnasiet in Gothenburg Sweden, Maria Nars, came up with the recursion relation in this example. I wanted to find a general formula and from that this article grew.} is, thus,  $1,\,2,\,3,\,3,\,6,\,18,\,54,\,324,\,...$\\

Let $b_n=\ln a_n$. Then $(b_n)_{n=0}^{\infty}$ is a number sequence with a linear recursion relation $b_{n+3}=b_{n+2}+b_n$ and can be treated like in ex 2 with initial values $\beta=\ln \alpha$, so $0,\,\ln 2,\,\ln 3$ and $N=3$. We will find $a_{16}$.

\[ b_n=\ln 2\cdot\kappa_n+(\ln 3 - \ln 2)\cdot\kappa_{n-1}=(\kappa_n-\kappa_{n-1})\ln 2+\kappa_{n-1}\ln 3\ , \]
\[ a_n=2^{\kappa_n-\kappa_{n-1}}\cdot 3^{\kappa_{n-1}} \]

 so we need $\kappa_{16}$ and $\kappa_{15}$, which have already been found in ex 2, $\kappa_{16}=189$ and $\kappa_{15}=129$. Thus,

\[ a_{16} = 2^{\kappa_{16}-\kappa_{15}}\cdot 3^{\kappa_{15}}=2^{189-129}\cdot 3^{129} \]
\[ =2^{60}\cdot 3^{129} \]
\[ =40779472028876430259264292468803306803871352789421825624677506478583962620919808 \]
(80 digits).\\

\vsp
\textbf{Ex 4:} $N=4\mbox{ and } \alpha=0,\,1,\,2,\,3$\ . \\

\[ a_n=\kappa_n+\kappa_{n-1}+\kappa_{n-2} \]

\vsp
\textbf{Ex 5:} $N=4\mbox{ and } \alpha=1,\,1,\,1,\,1$\ . \\

\[ a_n=\kappa_{n+1} \]

\vsp
\textbf{Ex 6:} $N=5\mbox{ and } \alpha=0,\,1,\,2,\,1,\,1$\ . \\

\[ a_n=\kappa_n+\kappa_{n-1}-\kappa_{n-2} \]

\section{Formulas for $D_n\mbox{ and } a_n$}

Partial fractions expansion of $f(x)$ (the zeroes $x_k$ of the denominator are simple):

\[ \frac{1}{1-x-x^N}=\sum_{k=1}^N\frac{r_k}{x-x_k} \]

where

\[ r_k=-\frac{1}{1+N {x_k}^{N-1}} \]

so, according to \eqref{fx},

\begin{eqnarray*}
  f(x) &=& \sum_{n=0}^{N-1}(\alpha_n-\alpha_{n-1})x^n\cdot\sum_{k=1}^N\frac{r_k}{x-x_k} \\
       &=& \sum_{n=0}^{N-1}(\alpha_n-\alpha_{n-1})x^n\cdot\sum_{k=1}^N\frac{-x_k}{(N-(N-1)x_k)(x-x_k)}\ \ \ .
\end{eqnarray*}

Since

\[ \dis\frac{\dop^m}{\dop x^m}x^{\ell}\Big |_{x=0} =m!\delta_{m\ell} \]

and

\[ \dis\frac{\dop^{m}}{\dop x^{m}}\frac{1}{x-x_k}\Big|_{x=0}=\frac{(-1)^m\cdot m!}{(-x_k)^{m+1}}=-\frac{m!}{{x_k}^{m+1}}\ , \]

\vsp
by Leibniz' formula for the $n$th derivative of a product

\begin{eqnarray*}
  \frac{\dop^n}{\dop x^n}\frac{x^{\ell}}{x-x_k}\Big |_{x=0}
   &=& \sum_{j=0}^{n} \binom{n}{j}\dis\frac{\dop^j}{\dop x^j}x^{\ell}\cdot\dis\frac{\dop^{n-j}}{\dop x^{n-j}}(x-x_k)^{-1}\Big |_{x=0} \\
   &=& -\binom{n}{\ell}\cdot \ell !\cdot (n-\ell)! {x_k}^{-1-n+\ell} \\
   &=& -n!\cdot {x_k}^{-1-n+\ell}\ \ \ .
\end{eqnarray*}

Since

\[ f(x) = \sum_{\ell=0}^{N-1}\sum_{k=1}^N \frac{(\alpha_{\ell}-\alpha_{\ell-1})(-x_k)}{N-(N-1)x_k}\cdot\frac{x^{\ell}}{x-x_k} \]

we get

\[ f^{(n)}(0)=\sum_{\ell=0}^{N-1}\sum_{k=1}^N \frac{(\alpha_{\ell}-\alpha_{\ell-1})x_k}{N-(N-1)x_k}\cdot n!\cdot {x_k}^{-1-n+\ell}\ . \]

According to \eqref{dnl}

\begin{eqnarray*}
  C_{n\ell} &=& \frac{1}{n!}\cdot\frac{\dop^n}{\dop x^n}\sum_{k=1}^{N}\frac{r_k x^{\ell}}{x-x_k} \Big |_{x=0} \\
   &=& \frac{1}{n!}\cdot\sum_{k=1}^{N}r_k(-n!){x_k}^{-1-n+\ell} \\
   &=& \sum_{k=1}^{N}\frac{1}{1+N {x_k}^{N-1}}\cdot{x_k}^{-1-n+\ell} \\
   &=& \sum_{k=1}^{N}\frac{{x_k}^{\ell-n}}{x_k+N{x_k}^N} \\
   &=& \sum_{k=1}^{N}\frac{{x_k}^{\ell}}{[N-(N-1)x_k]{x_k}^n}\ .
\end{eqnarray*}

Obviously $C_{n\ell}$ is a function of $n-\ell$ and can be written $C_{n\ell}=D_{n-\ell}$ where

\begin{equation}\label{Dn}
  D_n=\sum_{k=1}^{N}\frac{1}{[N-(N-1)x_k]{x_k}^n}
\end{equation}

and by \eqref{an1}, \eqref{dnl}

\begin{equation}\label{anf}
  a_n=\mathlarger{\mathlarger{\sum}_{k=1}^{N}\frac{\alpha_0+\dis\sum_{\ell=1}^{N-1}(\alpha_{\ell}-\alpha_{\ell-1}){x_k}^{\ell}}{[N-(N-1)x_k]{x_k}^n}}\ .
\end{equation}

Hence we have obtained explicit formulas for $D_n$ and $a_n$.

\section{Explicit function for $N=2$, Fibonacci again}

The polynomial $x^2+x-1$ has the zeroes $x_1=(-1+\sqrt{5})/2$, $x_2=(-1-\sqrt{5})/2$ and $x_1+x_2=-1$, $x_1\cdot x_2=-1$. $\alpha_0=0,\,\alpha_1=1$ and \eqref{anf} give\\

\begin{eqnarray*}
  F_n &=& \sum_{k=1}^2\,\frac{x_k}{(2-x_k){x_k}^n} \\
  \ \\
   &=& \dis\frac{1}{5(-1)^n}\left[{x_1}^n+{x_2}^n+2({x_1}^n x_2+x_1 {x_2}^n)\right] \\
   \ \\
   &=& \dis\frac{(-1)^n}{\sqrt{5}}\left({x_2}^n-{x_1}^n\right)
\end{eqnarray*}

after which substitution of $x_1$ och $x_2$ gives (Binet's formula)

\[ F_n=\dis\frac{1}{\sqrt{5}}\left\{\left(\dis\frac{1+\sqrt{5}}{2}\right)^n
-\left(\dis\frac{1-\sqrt{5}}{2}\right)^n\right\}\ . \]

\section{Explicit function for $N=3$}

Let

\[ p_n\equiv (x_1\,x_2)^n+(x_2\,x_3)^n+(x_3\,x_1)^n\ . \]

Then

\[ p_n=2(-\sqrt{\delta})^n\cdot T_n\left(\frac{\delta\sqrt{\delta}}{2}\right)+(\delta^2+1)^n \]

 where $T_n$, for $n=0,\,1,\,\dots$ denotes a Chebyshev polynomial\footnote{A Chebyshev polynomial $T_n$ is defined by $T_n(\cos\varphi)=\cos(n\varphi)$ or $T_n(x)=\cos(n\arccos x)$. The polynomial can be determined by means of trigonometric formulas for $\cos(n\varphi)$. For $n=0-5$ the polynomials are $T_0(x)=1,\,T_1(x)=x,\,T_2(x)=2x^2-1,\,T_3(x)=4x^3-3x,\,T_4(x)=8x^4-8x^2+1,\,T_5(x)=16x^5-20x^3+5x$. From the definition it follows that $|T_n(x)|\leq 1$ for all $x$.}

\vsp
\textbf{Proof for this:}\\

The zeroes of the polynomials $x^3+x-1$ satisfy

\begin{equation*}
  \left\{\begin{array}{rcl}
  x_1+x_2+x_3 &=& 0\ , \\
  \ \\
  x_1 x_2+x_2 x_3+x_3 x_1 &=& 1\ , \\
  \ \\
  x_1 x_2 x_3 &=& 1\ .
  \end{array}\right.
\end{equation*}

The real zero is denoted $x_1=\delta$ and the non-real ones can be written $x_2=\rho e^{i\varphi}$ and $x_3=\rho e^{-i\varphi}$. Since $x_1 x_2 x_3 = 1$, $\rho=1/\sqrt{\delta}$ and

\begin{eqnarray*}
  p_n &=& \left(\delta\cdot\rho e^{i\varphi}\right)^n + \rho^{2n}+\left(\delta\cdot\rho e^{-i\varphi}\right)^n \\
   \ \\
   &=& 2\delta^n\rho^n\cos(n\varphi)+\rho^{2n} \\
   \ \\
   &=& 2(\sqrt{\delta})^n\cos(n\varphi)+\delta^{-n} \\
   \ \\
   &=& 2(\sqrt{\delta})^n T_n(\cos\varphi)+(\delta^2+1)^n
\end{eqnarray*}

\vsp
 where $\cos\varphi=-\delta/(2\rho)=-\delta\sqrt{\delta}/2$. Since the Chebyshev polynomial $T_n$ has parity $(-1)^n$ we get

\[ p_n=2(-\sqrt{\delta})^n T_n(\delta\sqrt{\delta}/2)+(\delta^2+1)^n\ . \ \ \ \ \ \square \]

\vsp
We will now express $a_n$ in $p_n$ and start with $D_n$. $N=3$ in \eqref{Dn} gives\\

\[ D_n = \sum_{k=1}^3\,\frac{1}{(3-2x_k){x_k}^n}\ . \]

\vsp
At first we simplify the least common denominator,

\begin{eqnarray*}
   \prod_{k=1}^3\,(3-2x_k){x_k}^n &=& (3-2x_1)(3-2x_2)(3-2x_3)(x_1 x_2 x_3)^n \\
    \ \\
   &=& 27-18(x_3+x_1+x_2)+12(x_1 x_3+x_2 x_3+x_1 x_2)-8x_1 x_2 x_3\\
    \ \\
   &=& 27-18\cdot 0+12\cdot 1-8\cdot 1 = 31
\end{eqnarray*}

and then we get

\begin{eqnarray*}
  D_n &=& \dis\frac{1}{31}\left\{(3-2x_1){x_1}^n(3-2x_2){x_2}^n+(3-2x_2){x_2}^n(3-2x_3){x_3}^n+(3-2x_3){x_3}^n(3-2x_1){x_1}^n \right\} \\
  \ \\
   &=& \dis\frac{1}{31}\left\{(9-6(x_1+x_2)+4x_1 x_2)(x_1 x_2)^n+(9-6(x_2+x_3)+4x_2 x_3)(x_2 x_3)^n\right.\\
   \ \\
   & & \left.+(9-6(x_3+x_1)+4x_3 x_1)(x_3 x_1)^n\right\} \\
   \ \\
   &=& \dis\frac{1}{31}(9p_n+4p_{n+1}-6q_n)
\end{eqnarray*}

where

\[ q_n={x_1}^{n+1}{x_2}^n+{x_1}^n{x_2}^{n+1}+{x_2}^{n+1}{x_3}^n+{x_2}^n{x_3}^{n+1}
+{x_3}^{n+1}{x_1}^n+{x_3}^n{x_1}^{n+1}\ . \]

\vsp
Next we express $q_n$ in $p_n$. Since $x_1 x_2+x_2 x_3+x_3 x_1=1$,

\begin{eqnarray*}
  p_n &=& (x_1 x_2+x_2 x_3+x_3 x_1)\left[(x_1\,x_2)^n+(x_2\,x_3)^n+(x_3\,x_1)^n\right] \\
    \ \\
  &=& p_{n+1}+q_{n-1}\ .
\end{eqnarray*}

\vsp
Since this holds for all $n$ it can be written

\[ q_n=p_{n+1}-p_{n+2} \]

whereupon $D_n$ can be written

\[ D_n=\dis\frac{1}{31}(9p_n-2p_{n+1}+6p_{n+2}) \]

for $N=3$.\\

In the case $\alpha_0=0,\,\alpha_1=1,\,\alpha_2=2$ we get, by means of equation \eqref{an},

\begin{eqnarray*}
  a_n &=& D_{n-1}+D_{n-2} \\
  \ \\
  &=& \dis\frac{1}{31}(6p_{n+1}+4p_n+7p_{n-1}+9p_{n-2})\ .
\end{eqnarray*}

For example,

\begin{equation*}
  \begin{array}{rcl}
  a_3 &=& \dis\frac{1}{31}(6p_4+4p_3+7p_2+9p_1) \\
  \ \\

   &=& \dis\frac{1}{31}\{12\delta^2(8\cdot\frac{\delta^6}{16}-8\cdot\frac{\delta^3}{4}+1)
   -8\delta\sqrt{\delta}(4\cdot\frac{\delta^4\sqrt{\delta}}{8}-3\cdot\frac{\delta\sqrt{\delta}}{2})
   +14\delta(2\cdot\frac{\delta^3}{4}-1)\\
   \ \\
   & &-\,18\sqrt{\delta}\cdot\frac{\delta\sqrt{\delta}}{2}
   +6(\delta^2+1)^4+4(\delta^2+1)^3+7(\delta^2+1)^2+9(\delta^2+1)\} \\
   \ \\

   &=& \dots= \\
   \ \\
   &=& \dis\frac{1}{31}\{12\delta^2-24(1-\delta)+12(\delta-\delta^2)+24-48\delta+24\delta^2-24\delta^2
   +24(1-\delta) \\
   \ \\
   & & +\,62\delta-62\delta^2+12-12\delta+62\delta^2-14\delta+26\}=\dots =62/31=2\ .
     \end{array}
\end{equation*}

Obviously not any particularly efficient method for the calculation of particular terms $a_n$ of a number sequence.

\section{Estimate of $p_n$ for $N=3$ and one more formula for $a_n$}

From the above,

\[ p_n=2(-\sqrt{\delta})^n\cdot T_n\left(\frac{\delta\sqrt{\delta}}{2}\right)+(\delta^2+1)^n \]

where $\delta$ is the real zero of the polynomial $x^3+x-1$,

\[ \delta=\sqrt[3]{\dis\frac{\sqrt{93}}{18}+\frac{1}{2}}-\sqrt[3]{\dis\frac{\sqrt{93}}{18}-\frac{1}{2}} \]

and

\[ (\sqrt{\delta})^n<1/8\mbox{ if } n>10\ . \]

 Since $|T_n(x)|\leq 1$ because $T_n(x)$ is a cosine value, then $|p_n-(\delta^2+1)^n|<1/4$ from which follows that

\[ p_n=\lfloor(\delta^2+1)^n+1/2\rfloor \]

and thus $a_n$ can be written exactly

\[ a_n=\dis\frac{1}{31}\left\{6\lfloor(\delta^2+1)^{n+1}+1/2\rfloor+4\lfloor(\delta^2+1)^n+1/2\rfloor
+7\lfloor(\delta^2+1)^{n-1}+1/2\rfloor+9\lfloor(\delta^2+1)^{n-2}+1/2\rfloor\right\} \]

for $n>10$ or, alternatively,\\

\[ a_n=\dis\frac{1}{31}\left\{6\lfloor\epsilon^{n+1}+1/2\rfloor+4\lfloor\epsilon^n+1/2\rfloor
+7\lfloor\epsilon^{n-1}+1/2\rfloor+9\lfloor\epsilon^{n-2}+1/2\rfloor\right\} \]

where

\[ \epsilon=\dis\frac{1}{\delta}=\sqrt[3]{\dis\frac{29}{54}+\dis\frac{\sqrt{93}}{18}}
+\sqrt[3]{\dis\frac{29}{54}-\dis\frac{\sqrt{93}}{18}}+\dis\frac{1}{3}\ . \]

\vsp
Example:

\begin{eqnarray*}
  a_{16} &=& \dis\frac{1}{31}\left\{6\lfloor\epsilon^{17}+1/2\rfloor+4\lfloor\epsilon^{16}+1/2\rfloor
+7\lfloor\epsilon^{15}+1/2\rfloor+9\lfloor\epsilon^{14}+1/2\rfloor\right\} \\
\ \\
   &=& \dis\frac{1}{31}(6\cdot 664+4\cdot 453+7\cdot 309+9\cdot 211)=318\ .
\end{eqnarray*}

Corresponding formula for the sequence $(\kappa_n)_{n=0}^{\infty}$:

\begin{eqnarray*}
  \kappa_n &=& D_{n-1}=\dis\frac{1}{31}(9p_{n-1}-2p_n+6p_{n+1}) \\
  \ \\
   &=& \dis\frac{1}{31}\left\{6\lfloor\epsilon^{n+1}+1/2\rfloor-2\lfloor\epsilon^n+1/2\rfloor
   +9\lfloor\epsilon^{n-1}+1/2\rfloor\right\}\, .
\end{eqnarray*}

Similarly, but essentially simpler, the Fibonacci numbers can be expressed by means of the floor function $\lfloor\:\rfloor$:

\[ F_n=\left\lfloor\frac{1}{\sqrt{5}}\left(\frac{1+\sqrt{5}}{2}\right)^n+\dis\frac{1}{2}\right\rfloor \]

(for all $n$).\\

Also see \cite{graham}, p 300.

\section{General linear recursion relation}

$a_{n+N}=\sum_{u=1}^{N} A_u a_{n+N-u}$ and the initial values $(\alpha_n)_{n=0}^{N-1}$.\\

Generating function:

\begin{eqnarray*}
  f(x) &=& \sum_{n=0}^{\infty}a_n x^n=x^N\sum_{n=0}^{\infty}a_{n+N} x^n+\sum_{n=0}^{N-1}\alpha_n x^n \\
   \ \\
   &=& x^N\sum_{n=0}^{\infty}\sum_{u=1}^{N} A_u a_{n+N-u} x^n+\sum_{n=0}^{N-1}\alpha_n x^n \\
   \ \\
   &=& \sum_{u=1}^{N-1}A_u x^u \sum_{n=0}^{\infty} a_{n+N-u} x^{n+N-u}+A_N x^N \sum_{n=0}^{\infty}a_n x^n+\sum_{n=0}^{N-1}\alpha_n x^n\\
   \ \\
   &=& \sum_{u=1}^{N}A_u x^u\cdot f(x)-\sum_{u=1}^{N-1}\sum_{n=0}^{N-u-1}A_u\alpha_n x^{u+n}+\sum_{n=0}^{N-1}a_n x^n
 \end{eqnarray*}\\

giving

\[ \left(1-\sum_{n=1}^N A_n x^n\right)\cdot f(x)
=\sum_{n=0}^{N-1}\alpha_n x^n - \sum_{u=1}^{N-1}\sum_{n=0}^{N-u-1}A_n \alpha_n x^{u+n}\ . \]\\

Example: $N=3$\\

\begin{eqnarray*}
  \sum_{u=1}^{N-1}\sum_{n=0}^{N-u-1}A_u\alpha_n x^{u+n}
  &=& A_1\alpha_0 x+A_1\alpha_1 x^2+A_2\alpha_0 x^2
\end{eqnarray*}

from which

\begin{equation*}
  f(x)=\frac{x^2(\alpha_2-A_1\alpha_1-A_2\alpha_0)+x(\alpha_1-A_1\alpha_0)+\alpha_0}{1-\sum_{n=1}^{3}A_n x^n}\ .
\end{equation*}

\section{Conclusion}

Given a relation $a_n=a_{n-1}+a_{n-N}$ and $N$ initial values you can, obviously, calculate any value using the recursion relation repeatedly. But a closed formula would directly give the desired value. This, however, proved rather intricate when $N>2$ since the zeroes of a polynomial of degree $N$ will occur in the formula. A formula derived for $N=3$ is relatively simple but yet tedious to use for large index values. Using, however, a suitable computer program, such as DERIVE, it will work well. Apart from this a closed formula is, of course, interesting in itself. Further, a method has been deduced making it possible to calculate any desired value $a_n$ by associating the value to another number sequence $(\kappa_n)_{n=0}^{\infty}$, which is easier to calculate. This method, which I regard as the most interesting thing of this article, is an inbetween of directly iterating the recursion formula and using an expression in closed form.\\

Another method of finding an explicit formula for the terms in a number sequence, besides using a generating function, is to determine the zeroes of the characteristic polynomial for the recursion relation, analogously with a method of solving linear, ordinary differential equations. For $N=2$ the amount of work as well as the formula obtained will be reasonably equal, but for $N=3$ essentially worse as far as I have investigated. It would be of some interest to compare the methods in general. The algorithm using the sequence $(\kappa)$ works principally the same way independently of $N$. Finally, it might be interesting to investigate the influence of different $N$ on the zeroes of the polynomial $x^N+x-1$ or the characteristic polynomial $x^N-x^{N-1}-1$, which also has degree $N$.\footnote{These polynomial have a simple connection: If one polynomial has zeroes $x_k$, then the other has zeroes $1/x_k$. They are essentially reciprocal polynomials of each other.}

\renewcommand\appendixname{Appendix}

\begin{appendices}

\section*{Appendix\\ \\The initial values for ${\kappa}$}

Proof that $\kappa_n=1$ for $0<n<N$ ($N\geq 2$).\\

Notation: $\kappa_n\equiv \kappa^N_{n,N-1}$.\\

\begin{eqnarray*}
  (\delta^{N-1}+1)^n &=& \sum_{\nu=0}^{n}\binom{n}{\nu}\delta^{\nu N-\nu} \\
  \ \\
   &=& \sum_{\nu=0}^{n}\binom{n}{\nu}(1-\delta)^{\nu-1}\delta^{N-\nu} \\
   \ \\
   &=& 1+\sum_{\nu=1}^{n}
   \sum_{\mu=0}^{\nu-1}\binom{n}{\nu}\binom{\nu-1}{\mu}(-1)^{\mu}\delta^{N-\nu+\mu}
\end{eqnarray*}

where

\[ 1\leqslant N-n\leqslant N-\nu\leqslant N-\nu+\mu\leqslant N-\nu+\nu-1=N-1 \]

whence

\[ 1\leqslant N-\nu+\mu\leqslant N-1\ . \]\\

The coefficient of $\delta^{N-1}$, i e $\kappa_n=\kappa^N_{n,N-1}$ is obtained for $\mu-\nu=-1$, i e

\begin{eqnarray*}
   \kappa_n&=& \sum_{\nu=1}^{n}\binom{n}{\nu}\binom{\nu-1}{\nu-1}(-1)^{\nu-1} \\
   \ \\
   &=& -\sum_{\nu=0}^{n}\binom{n}{\nu}(-1)^{\nu}+\binom{n}{0}(-1)^0 \\
   \ \\
   &=& -(1-1)^n+1=1\ .
\end{eqnarray*}

\end{appendices}


\end{document}